\documentclass[12pt,a4paper,reqno]{amsart}
\usepackage{amsmath}
\usepackage[latin1]{inputenc}
\usepackage{amsfonts}
\usepackage{amssymb}
\numberwithin{equation}{section}

\def\R{\mathbb{R}}
\def\E{\mathbb{E}}
\def\dint{\displaystyle\int}

\def\P{\mathbb{P}}

\newenvironment{proof1}[1][Proof]{\emph{#1.} }{\rule{0.5em}{0.5em}}

\newenvironment{proof2}[1][Proof]{\emph{#1.} }{\rule{0.5em}{0.5em}}

\theoremstyle{plain}

\newtheorem{thm}[subsection]{Theorem}

\newtheorem{prop}[subsection]{Proposition}

\newtheorem{lem}[subsection]{Lemme}
\newtheorem{rem}[subsection]{Remarque}

\theoremstyle{definition}

\theoremstyle{remark}

\renewcommand{\leq}{\leqslant}
\renewcommand{\geq}{\geqslant}

\newsavebox{\proofbox}
\savebox{\proofbox}{\begin{picture}(7,7)%
  \put(0,0){\framebox(7,7){}}\end{picture}}

\def\E{\mathbb{E}}

\def\R{\mathbb{R}}

\def\P{\mathbb{P}}

\begin{document}
\title[]{Unicité trajectorielle des équations différentielles stochastiques avec temps local et temps de séjour au bord}
\address[R. Belfadli and Y. Ouknine]{Department of Mathematics, Faculty of Sciences Semlalia\\
Cadi Ayyad University, B.P. 2390 Marrakesh, Morocco.}
\email[R. Belfadli]{r.belfadli@ucam.ac.ma}
\author{R. Belfadli}
\author{Y. Ouknine}
\address[Y. Ouknine]
{Hassan II Academy of Sciences and
Technology.}
\email[Y. Ouknine]{ouknine@ucam.ac.ma}
\thanks{This research is
supported by the Hassan II Academy of Sciences and Technology.}
\date{\today}
\begin{abstract}
Nous étudions l'unicité trajectorielle des solutions d'une classe d'équations différentielles stochastiques avec temps local et temps de séjour au bord. Nous utilisons le problème des martingales associé pour montrer qu'il y a unicité en loi, puis nous établissons  que le supremum de deux solutions est encore une solution.
\end{abstract}

\maketitle
\section{Introduction}
Nous nous intéressons dans ce papier à l'étude de l'unicité trajectorielle des solutions d'équations différentielles stochastiques (EDSs) avec temps local et temps de séjour au bord. Plus précisement, on considère sur un espace probabilisé  $(\Omega,\mathcal{F},\mathbb{P})$ les EDSs du type:
\begin{align}\label{W-G}
 \left\{
\begin{aligned}
&X_{t} = x+ \int_{0}^{t}1\!_{\{X_s \neq 0\}}dB_s + \int_{0}^{t}\alpha(s) dL^{0}_s(X)\\
&\int_{0}^{t}1\!_{\{X_s= 0\}}ds =  \int_{0}^{t}\rho(h_s)dL^{0}_{s}(X),
\end{aligned}
\right.
\end{align}
où dans cette équation $B$ désigne un mouvement brownien linéaire issu de $0$, $L^{0}_{t}(X)$ est le temps local au point zéro de la semimartingale inconnue $X$, $\alpha(\cdot)$ est une fonction borélienne positive,  $\rho(\cdot)$ est une fonction lipschitzienne strictement positive sur $[0,1]$ et la fonction inconnue $h_t$ est définie par $h_t:=\P_{x}(X_t= 0)$. \vspace{.2cm}\\
Dans cette équation, il ya en fait deux inconnues: le processus aléatoire à valeurs réelles $(X_t, t \geq 0)$ et la fonction $h_t$.\vspace{.2cm}\\
Lorsque la fonction $\alpha(\cdot)$ est une constante inférieur à $1/2$, (\ref{W-G}) se réduit à une EDS introduite par S. Weinryb dans \cite{S-W}. Elle posséde, dans ce cas, une unique solution trajectorielle obtenue comme le ``processus non linéaire" associé à un système de
particules en intéraction dont le temps de séjours de chacune dépend du nombre moyen de particules au bord par l'intermidiaire d'une fonction
décroissante $\rho$ (voir, \cite{S-W} pour plus de détails). Indiquons également qu'on retrouve ces EDSs dans S. Watanabe \cite{ZW} lors de
l'étude de quelques exemples explicites d'EDSs.
 \vspace{.15cm}\\
Le type d'EDSs (\ref{W-G}), que nous considérons ici, est légèrement plus général que celui considéré par S. Weinryb dans \cite{S-W}.
Notre objectif est de montrer l'unicité trajectorielle des solutions à ces EDSs. Pour cela, nous allons utiliser
de façon essentielle une technique introduite par A. Y. Veretennikov \cite{ver} et utilisée par E. Perkins \cite{per}, J. F. Le Gall \cite{Legall} et S. Weinryb
(\cite{sw}, \cite{S-W}), qui consiste à démontrer l'unicité en loi pour ces solutions, et que le supremum de deux solutions est encore une solution.\vspace{.2cm}\\
\noindent Nous procédons donc comme suit: Dans la Section \ref{P-M1}, on formule le problème des martingales associé au système (\ref{W-G}), puis on montre l'unicité de la solution à ce problème.
Ensuite, dans la Section \ref{P-M2}, nous utilisons le temps local pour montrer que le maximum et le minimum de deux solutions sont encore des
solutions.
\vspace{.2cm}\\
\section{Formulation du problème martingale non linéaire associé}\label{P-M1}
Sur l'espace canonique $C(\R)$ on s'intéresse aux solutions $\P_{x}$ du problème de martingale $\mathcal{P}(x)$:
\begin{enumerate}
  \item pour toute fonction $f \in C^{2}(\R^{*})\cap C(\R)$ possédant des dérivées à droites et à gauche en $0$, bornée ainsi que ses dérivées
  \[f(X_t)-f(X_0)-\frac{1}{2}\int_{0}^{t}f^{\prime\prime}(X_s)1\!_{\{X_s \neq 0\}}ds -\int_{0}^{t}[f^{\prime}(0^{+})\alpha(s)+ \frac{f^{\prime}(0^{+})-f^{\prime}(0^{-})}{2}]dL_s^{0}(X)\] est une martingale. Ici, $L^{0}_{t}(X)$ désigne le temps local au sens de Tanaka de la semimartingale $X$;
  \item $\P_x(X_0=x)=1$ et $\int_{0}^{t}1\!_{\{X_s= 0\}}ds = \int_{0}^{t}\rho(h_s)dL^{0}_{s}(X).$
\end{enumerate}
\vspace{.3cm}
\textbf{Unicité faible de la solution $\P_{x}$ de $\mathcal{P}(x)$:}\\
  \ Nous montrerons tout d'abord que la loi $\mu_{t}(x, dy):=\P_{x}(X_t\in dy)$ est absolument continue par rapport à
la mesure de Lebesgue $dy$ sur $\R^{*}$.\\
\begin{lem}\label{lem1}
Pour tout $t\geq 0$, on a
\begin{eqnarray}\label{abso}
  \lim_{\epsilon\rightarrow 0}\int_{0}^{t}\frac{\mu_{t}(x, A_{\epsilon})}{2\epsilon}ds=\int_{0}^{t}\frac{1-\alpha(s)}{\rho(h_s)}h(s)\ ds,
\end{eqnarray}
où $A_{\epsilon}:= ]-\epsilon,\epsilon[\setminus\{0\}$ et $h(t):= \P_{x}(X_t=0)$.
En particulier, $\mu_{t}(x, dy)$ est absolument continue par rapport à la mesure de Lebesgue sur $\R^{*}$ et sa densité de Radon-Nikodym $p_{t}(x,y)$ vérifie
\begin{eqnarray}\label{formule1}
  \frac{p_{t}(x,0^{+})+ p_{t}(x,0^{-})}{2}=\frac{1-\alpha(t)}{\rho(h_t)}h(t)
\end{eqnarray}
\end{lem}
\noindent\begin{proof1}[Démonstration du Lemme \ref{lem1}]
D'aprés le point $(i)$, on a pour toute fonction $f$ régulière
\[\E\left(\int_{0}^{t} f^{\prime\prime}(X_s)1\!_{\{X_s \neq 0\}}ds \right)= \E\left(\int_{0}^{t} f^{\prime\prime}(X_s)d\langle X,X \rangle  _{s}\right)\]
ce qui implique $d\langle X,X \rangle  _{s}=1\!_{\{X_s \neq 0\}}ds$, et par conséquent
\begin{align}\label{equa1}
 \frac{1}{2}(L^{0}_{t}+L^{0^{-}}_{t}) =\lim_{\epsilon\rightarrow 0}\frac{1}{2\epsilon}\int_{0}^{t} 1\!_{\{0\leq|X_s|\leq \epsilon\}} d\langle X,X \rangle  _{s}=\lim_{\epsilon\rightarrow 0}\frac{1}{2\epsilon}\int_{0}^{t} 1\!_{\{0<|X_s|\leq \epsilon\}} ds.
\end{align}
D'autre part,
\[\frac{1}{2}(L^{0}_{t}-L^{0^{-}}_{t})=\int_{0}^{t} 1\!_{\{X_s=0\}}dV_s,\]
où $V_t= \int_{0}^{t} \alpha(s)dL^{0}_{s}(X)$ désigne la partie à variation finie de la semimartingale $X$. Ce qui donne
$L^{0^{-}}_{t}=\int_{0}^{t}(1-2\alpha(s))dL^{0}_{s}$. Et par suite, compte tenu du point $(ii)$, on obtient:
\begin{align}\label{equa2}
 \frac{1}{2}(L^{0}_{t}+L^{0^{-}}_{t}) =\int_{0}^{t}(1-\alpha(s))dL^{0}_{s}=\int_{0}^{t}\frac{1-\alpha(s)}{\rho(h_s)}1\!_{\{X_s=0\}} ds.
\end{align}
D'où, l'égalité (\ref{abso}) s'obtient en prenant l'espérance des seconds membres des deux équations (\ref{equa1}) et (\ref{equa2}).
\end{proof1}
\begin{lem}\label{lem22}
La fonction $h(t):= \P_{x}(X_t=0)$ satisfait l'équation intégro-différentielle
\begin{eqnarray}\label{int-diff}
  \frac{1-\alpha(t)}{\rho(h_t)}\!\ h(t) = \frac{1}{\sqrt{2\pi t}}\!\ e^{-\frac{x^{2}}{2 t}} + \frac{1}{2\sqrt{2\pi}} \dint_{0}^{t} \frac{h(t-s)-h(t)}{s^{3/2}}\!\ ds - \frac{h(t)}{\sqrt{2\pi t}}
\end{eqnarray}
\end{lem}
\noindent\begin{proof1}[Démonstration du Lemme \ref{lem22}]Posons, pour $t\geq0$ et $\lambda\in \R$,
\begin{eqnarray*}
f(\lambda,t)=\E_x(e^{i\lambda X_t}1_{\{X_t\neq 0\}})
\end{eqnarray*}

Par application de la formule d'Itô et en tenant compte que la partie à variation finie $V$ du processus $X$ est donnée par
$dV_t=\alpha(t)dL^{0}_{t}(X)=\frac{\alpha(t)}{\rho(h_s)}1_{\{X_t=0\}}dt$, on a:
\begin{eqnarray*}
  f(\lambda,t)+ h(t)= e^{i\lambda x}+ i\lambda\  \E_{x}\left(\int_{0}^{t}e^{i\lambda X_s}\frac{\alpha(s)}{\rho(h_s)}\!\ 1_{\{X_s=0\}}ds\right)-
  \frac{\lambda^{2}}{2}\int_{0}^{t}f(\lambda,s)ds
\end{eqnarray*}
Soit encore,
\begin{eqnarray*}
f(\lambda,t)+ h(t)=e^{i\lambda x}+  \dint_{0}^{t}\left[\frac{i\lambda}{\rho(h_s)}\alpha(s)+ \frac{\lambda^{2}}{2}\right]h(s) ds-
  \frac{\lambda^{2}}{2}\dint_{0}^{t}(f(\lambda,s)+h(s))ds
\end{eqnarray*}
En considérant la fonction $h$ comme connue dans cette équation et aprés résolution, on obtient:
\begin{eqnarray*}
f(\lambda,t)+ h(t)&=e^{-\frac{\lambda^{2}}{2}t}\left[e^{i\lambda x}-\frac{\lambda^{2}}{2}\int_{0}^{t}ds\ e^{\frac{\lambda^{2}}{2}s}\int_{0}^{s}du\left(\frac{i\lambda}{\rho(h_s)}\!\ \alpha(u)+\frac{\lambda^{2}}{2}\right)h(u)\right]\\
&\  + \int_{0}^{t}dsh(s)\left(\frac{i\lambda}{\rho(h_s)}\alpha(s)+\frac{\lambda^{2}}{2}\right)
\end{eqnarray*}
Soit, aprés une intégration par partie,
\begin{eqnarray}\label{loi1}
f(\lambda,t)=e^{-\frac{\lambda^{2}}{2}t}\left[e^{i\lambda x}+\int_{0}^{t}ds\ e^{\frac{\lambda^{2}}{2}s}h(s)\left(\frac{i\lambda}{\rho(h_s)}\!\ \alpha(s)+\frac{\lambda^{2}}{2}\right)\right]- h(t).
\end{eqnarray}
D'autre part, $f(\lambda,t)$ est la transformée de Fourier de la fonction $p_{t}(x,y)$ et donc il s'ensuit par la formule d'inversion que
\[\frac{p_{t}(x,0^{+})+ p_{t}(x,0^{-})}{2}= \frac{1}{2\pi}\int f(\lambda, t)d\lambda\]
En injectant (\ref{loi1}) dans cette formule et en tenant compte de (\ref{formule1}), on obtient (\ref{int-diff}).
\end{proof1}
\begin{prop}\label{prop1}
 Lorsque la fonction $\alpha(t)\leq 1/2$, il ya unicité faible de la solution $\P_x$ de $\mathcal{P}(x)$.
\end{prop}
\begin{proof2}[Preuve de la Proposition \ref{prop1}]
Observant tout d'abord que, puisque $\alpha(\cdot)$ est bornée et que la fonction $g(t,x):=\frac{1-\alpha(t)}{\rho(x)}x$ définissant le second membre de (\ref{int-diff}) est lipshitzienne en $x$, l'argument de S. Weinryb \cite{S-W} s'adapte bien et on a l'unicité de la solution de l'équation intégro-différentielle (\ref{int-diff}). On en déduit alors, de l'expression (\ref{loi1}), que la fonction $f(\lambda, t)$ est uniquement déterminée. Autrement dit, on a l'unicité de la loi de $X_t$ sous $\P_{x}$ pour chaque $t>0$. Puis on passe à tout le processus en utilisant le caractère markovien de $X$.
\end{proof2}
\section{Unicité trajectorielle des solutions de (\ref{W-G})}\label{P-M2}
\begin{thm}
Si la fonction $\alpha(t)\leq \frac{1}{2}$, alors il ya unicité trajectorille des solutions de (\ref{W-G}).
\end{thm}
\begin{proof2}[Preuve du théorème]
D'aprés la première partie, il ya unicité en loi des solutions de (\ref{W-G}). Il suffit alors de démontrer que le supremum de deux solutions est encore une solution.\\
 Si $X^{1}$ et $X^{2}$ sont deux solutions de (\ref{W-G}), associées au même brownien $B$, alors en écrivant $X^{1}_{t}\vee X^{2}_{t}=(X_{t}^{1}-X_{t}^{2})^{+}+X_{t}^{2}$ et par application de la formule d'Itô-Tanaka, on ~a
\begin{eqnarray}\label{supdecomp}
d(X^{1}\vee X^{2})_{t}&=&1\!_{\{X^{1}_{t}\vee X^{2}_{t}\neq 0\}}dB_t + \alpha(t)\left(1\!_{\{X^{2}_{t}<0\}}dL^{0}_{t}(X^{1})
+1\!_{\{X^{1}_{t}\leq 0\}} dL^{0}_{t}(X^{2}) \right)\nonumber\\
&& \quad \quad +\frac{1}{2}dL^{0}_{t}(\Delta X),
\end{eqnarray}
où $L^{0}_{t}(\Delta X)$ désigne le temps local de la semimartingale $\Delta X=X^{1}-X^{2}$.

\noindent De plus,
\begin{eqnarray}\label{bordsup}
  1\!_{\{X^{1}_{t}\vee X^{2}_{t}=0\}}dt &=& \rho(h_s)\left(1\!_{\{X^{2}_{t}<0\}}dL^{0}_{t}(X^{1})
+1\!_{\{X^{1}_{t}\leq 0\}} dL^{0}_{t}(X^{2})\right).
\end{eqnarray}
Ce qui montre alors que le supremum $(X^{1}_{t}\vee X^{2}_{t}, t\geq 0)$ est également solution de (\ref{W-G}) grâce au résultat de la Proposition \ref{tempsup} ci-dessous et  par suite l'unicité trajectorielle pour l'équation $(\ref{W-G})$.
\end{proof2}\\
\begin{prop}\label{tempsup}Pour tout $t\geq 0$, on a
\begin{eqnarray}
L^{0}_{t}(\Delta X)=0 \quad \mbox{et}\quad dL^{0}_{t}(X^{1}\vee X^{2})=1\!_{\{X^{2}_{t}<0\}}dL^{0}_{t}(X^{1})
+1\!_{\{X^{1}_{t}\leq 0\}} dL^{0}_{t}(X^{2}).
\end{eqnarray}
\end{prop}
\vspace{.6cm}
Pour la preuve de cette proposition, nous avons besoin de deux lemmes.
\vspace{.3cm}\\
\noindent Introduisons, pour $0<\alpha<\frac{1}{2}$ et $i=1, 2$, les processus $(Z^{\alpha, i}_{t}:=X^{i}-2\alpha X^{i +}_{t}, t\geq 0)$ et
$(\Delta Z^{\alpha}_{t}:=Z^{\alpha, 1}_{t}-Z^{\alpha, 2}_{t}, t\geq 0)$.

\begin{lem}\label{1}
\begin{enumerate}
  \item [($\imath$)] La mesure $dL_{t}^{0}(X^{1}-X^{2})$ est absolument continu par rapport à la mesure de Lebesgue et à support inclu dans l'ensemble $\{t, X^{1}_{t}=X^{2}_{t}=0\}$;
  \item [($\imath\imath$)] Il existe une fonction mesurable positive $\varphi$ telle que, pour tout $t\geq 0$
  \begin{eqnarray}\label{repr}
   L_{t}^{0}(\Delta Z^{\alpha})=\int_{0}^{t}(1-2\alpha \varphi(s))dL_{s}^{0}(\Delta X)
  \end{eqnarray}
  \item [($\imath\imath\imath$)] Si $\varphi$ est la fonction définie en $(\imath\imath)$, alors on a pour tout $t\geq0$
  \begin{eqnarray}\label{templocalsup1}
    L^{0}_{t}(X^{1}\vee X^{2})= \int_{0}^{t}1\!_{\{X^{1}_{s}\leq 0\}} dL^{0}_{s}(X^{2})+\int_{0}^{t}1\!_{\{X^{2}_{s}<0\}}dL^{0}_{s}(X^{1})
 + \int_{0}^{t}\varphi(s)dL_{s}^{0}(\Delta X)
  \end{eqnarray}
\end{enumerate}
\end{lem}
\begin{rem}De l'expression (\ref{templocalsup1}), on en déduit que
\[\int_{0}^{t}\varphi(s)dL_{s}^{0}(\Delta X)=\frac{1}{2\alpha}[L_{t}^{0}(\Delta X)-L_{t}^{0}(\Delta Z^{\alpha})]\]
est indépendant de $\alpha$.

\end{rem}
\begin{proof1}[Preuve du Lemme \ref{1}]($\imath$) L'absolu continuité résulte de l'estimation
\[L^{0}_{t}(\Delta X)\leq \dint_{0}^{t}\left(1\!_{\{X^{1}_s= 0\}}+1\!_{\{X^{2}_s= 0\}}\right)\frac{ds}{\rho(h_s)}, \quad \forall t\geq 0. \]
Cette inégalité est une conséquence immédiate de $L^{0}_{t}(\Delta X)\leq L^{0}_{t}(X^{1})+L^{0}_{t}(X^{2})$ et du fait que $X^{1}$ et $X^{2}$ sont solutions de (\ref{W-G}). Pour montrer que le support de $dL^{0}_{t}(\Delta X)$ est inclu dans $\{t, X^{1}_{t}=X^{2}_{t}=0\}$, il suffit de remarquer que
$\int_{0}^{t} |X^{1}_{s}|dL^{0}_{s}(\Delta X)=0.$\vspace{.2cm}\\
%
($\imath\imath$) Par application du lemme $1$ de \cite{ouknak}, on a $L_{t}^{0}(\Delta Z^{\alpha})\leq L_{t}^{0}(\Delta X)$. Il en résulte que $dL_{t}^{0}(\Delta Z^{\alpha})$ est absolumlent continue par rapport à $dL_{t}^{0}(\Delta X)$ et par suite on a la représentation (\ref{repr}).
\vspace{.2cm}\\
($\imath\imath\imath$) D'aprés \cite{Revuz-Yor} (voir aussi \cite{OR}), on sait que
\begin{eqnarray*}
    L^{0}_{t}(X^{1}\vee X^{2})= \int_{0}^{t}1\!_{\{X^{1}_{s}\leq 0\}} dL^{0}_{s}(X^{2})+ L_{t}^{0}(X^{1+}-X^{2+})
  \end{eqnarray*}
 il suffit donc de montrer que
  \[L_{t}^{0}(X^{1+}-X^{2+})=\int_{0}^{t}1\!_{\{X^{2}_{s}< 0\}} dL^{0}_{s}(X^{1})+ \int_{0}^{t}\varphi(s)dL_{s}^{0}(\Delta X)\]
  où encore, compte tenu de (\ref{repr}), que
  \begin{eqnarray}\label{tt}
  L_{t}^{0}(X^{1+}-X^{2+})=\int_{0}^{t}1\!_{\{X^{2}_{s}< 0\}} dL^{0}_{s}(X^{1})+ \frac{1}{2\alpha}[L_{t}^{0}(\Delta X)-L_{t}^{0}(\Delta Z^{\alpha})].
  \end{eqnarray}
 En utilisant l'égalité
 \[2\alpha (X_{t}^{1+}-X_{t}^{2+})^{+}=(\Delta X_{t})^{+}-(\Delta Z_{t}^{\alpha})^{+}\]
on a, d'une part, par la formule d'Itô-Tanaka
\begin{eqnarray}\label{equality1}
2\alpha\hspace{.08cm} d(X_{t}^{1+}-X_{t}^{2+})^{+}=2\alpha 1\!_{\{X^{1+}_{t}>X^{2+}_{t}\}}d(X_{t}^{1+}-X_{t}^{2+})+ \alpha \hspace{.08cm}dL^{0}_{t}(X_{t}^{1+}-X_{t}^{2+})
\end{eqnarray}
et d'autre part,
  \begin{eqnarray*}
  2\alpha\hspace{.08cm} d(X_{t}^{1+}-X_{t}^{2+})^{+}=1\!_{\{\Delta X_{t}>0\}}d(\Delta X_{t})-1\!_{\{\Delta Z_{t}^{\alpha}>0\}}d(\Delta Z_{t}^{\alpha})+\frac{1}{2}d(L_{t}^{0}(\Delta X)-L_{t}^{0}(\Delta Z^{\alpha}))
  \end{eqnarray*}
  or $\Delta X_{t}\Delta Z^{\alpha}_{t}\geq0$, il vient
  \begin{eqnarray}\label{equality2}
  2\alpha\hspace{.08cm} d(X_{t}^{1+}-X_{t}^{2+})^{+}=2\alpha\hspace{.08cm}1\!_{\{\Delta X_{t}>0\}}d(X_{t}^{1+}-X_{t}^{2+})+\frac{1}{2}\hspace{.08cm}d(L_{t}^{0}(\Delta X)-L_{t}^{0}(\Delta Z^{\alpha})).
  \end{eqnarray}
 Ce qui conduit, par (\ref{equality1}), à
  \begin{eqnarray}\label{mm}
dL^{0}_{t}(X_{t}^{1+}-X_{t}^{2+})&=&2\left(1\!_{\{X^{1}_{t}>X^{2}_{t}\}}-1\!_{\{X^{1+}_{t}>X^{2+}_{t}\}}\right)\hspace{.08cm} d(X_{t}^{1+}-X_{t}^{2+})\nonumber\\
&&\quad \quad+\frac{1}{2\alpha}\hspace{.08cm}d(L_{t}^{0}(\Delta X)-L_{t}^{0}(\Delta Z^{\alpha})).
\end{eqnarray}
Par ailleurs, en appliquant la formule d'Itô-Tanaka, on montre aisément que
  \begin{eqnarray}\label{equality2}
  2\left(1\!_{\{X^{1}_{t}>X^{2}_{t}\}}-1\!_{\{X^{1+}_{t}>X^{2+}_{t}\}}\right)\hspace{.08cm} d(X_{t}^{1+}-X_{t}^{2+})=1\!_{\{X^{2}_{t}< 0\}} dL^{0}_{t}(X^{1})
  \end{eqnarray}
  et donc (\ref{tt}) s'obtient en reportant (\ref{equality2}) dans (\ref{mm}). Ceci complète la preuve du Lemme \ref{1}.
\end{proof1}
\begin{lem}\label{lem2}
\begin{enumerate}
  \item[($\imath$)] Pour tout $t\geq0$, on a
  \begin{eqnarray}\label{intphi}
   \int_{0}^{t}\varphi(s)dL_{s}^{0}(\Delta X)= 0
  \end{eqnarray}
  \item[($\imath\imath$)] Posons: $\quad$ $h_{1,2}(t)=\P_{x}(X^{1}_{t}\vee X_{t}^{2}=0) \quad \mbox{et} \quad \overline{h}_{1,2}(t)=\P_{x}(X^{1}_{t}\wedge X_{t}^{2}=0)$. Alors $h_{1,2}$ et $\overline{h}_{1,2}$ satisfont aux équations suivantes:

 \begin{eqnarray}\label{int-diff2}
  \frac{1-\alpha(t)}{\rho(h_t)}\!\ h_{1,2}(t)\hspace{-.1cm} +\hspace{-.1cm}  \E\psi^{\prime}(t)=\frac{1}{\sqrt{2\pi t}}\!\ e^{-\frac{x^{2}}{2 t}} +\hspace{-.1cm} \frac{1}{2\sqrt{2\pi}}\hspace{-.08cm} \dint_{0}^{t} \frac{h_{1,2}(t-s)-h_{1,2}(t)}{s^{3/2}}\!\ ds -\hspace{-.08cm} \frac{h_{1,2}(t)}{\sqrt{2\pi t}}
\end{eqnarray}
\begin{eqnarray}\label{int-diff3}
  \frac{1-\alpha(t)}{\rho(h_t)}\!\ \overline{h}_{1,2}(t)\hspace{-.1cm}+\hspace{-.1cm}\E\psi^{\prime}(t)=\frac{1}{\sqrt{2\pi t}}\!\ e^{-\frac{x^{2}}{2 t}} +\hspace{-.1cm} \frac{1}{2\sqrt{2\pi}}\hspace{-.08cm} \dint_{0}^{t} \frac{\overline{h}_{1,2}(t-s)-\overline{h}_{1,2}(t)}{s^{3/2}}\!\ ds -\hspace{-.08cm} \frac{\overline{h}_{1,2}(t)}{\sqrt{2\pi t}}
\end{eqnarray}
Ici, la fonction $\psi$ est définie par $\psi(t)=-\frac{1}{2}L_{t}^{0}(\Delta X)$ et $\psi^{\prime}$ désigne sa dérivée de Radon-Nikodym par rapport à la mesure de Lebesgue.
\end{enumerate}
\end{lem}
\begin{proof1}[Preuve du Lemme \ref{lem2}]
($\imath$) Puisque
\[d\langle \Delta Z^{\alpha},\Delta Z^{\alpha} \rangle  _{t}=\left[1\!_{\{X^{1}_t \neq 0\}}-1\!_{\{X^{2}_t \neq 0\}}-2\alpha(1\!_{\{X^{1}_t > 0\}}-1\!_{\{X^{2}_t > 0\}})\right]^{2}dt\]
on a
\begin{align*}
 L_{t}^{0}(\Delta Z^{\alpha})=\lim_{\epsilon\rightarrow 0}\frac{1}{2\epsilon}\int_{0}^{t} 1\!_{\{0\leq\Delta X_{s}-2\alpha\Delta Z^{\alpha}_{s}\leq \epsilon\}} \left[1\!_{\{X^{1}_s \neq 0\}}-1\!_{\{X^{2}_s \neq 0\}}-2\alpha(1\!_{\{X^{1}_s > 0\}}-1\!_{\{X^{2}_s > 0\}})\right]^{2}ds
\end{align*}
En développant le carré, il vient
\begin{align}\label{appro}
 L_{t}^{0}(\Delta Z^{\alpha})=\lim_{\epsilon\rightarrow 0}\left[I^{1}_{\epsilon}(t)+I^{2}_{\epsilon}(t)+I^{3}_{\epsilon}(t)\right],
\end{align}
où
\begin{align*}
 I^{1}_{\epsilon}(t)&=\frac{1}{\epsilon}\int_{0}^{t} 1\!_{\{0\leq\Delta X_{s}-2\alpha\Delta Z^{\alpha}_{s}\leq \epsilon\}} \left[1\!_{\{X^{2}_s = 0\}}-1\!_{\{X^{1}_s = 0\}}\right]^{2}ds;\\
 I^{2}_{\epsilon}(t)&=4\alpha^{2}\frac{1}{\epsilon}\int_{0}^{t} 1\!_{\{0\leq\Delta X_{s}-2\alpha\Delta Z^{\alpha}_{s}\leq \epsilon\}}\left[1\!_{\{X^{1}_s > 0\}}-1\!_{\{X^{2}_s > 0\}}\right]^{2} ds;
 \end{align*}
 et
 \begin{align*}
 I^{3}_{\epsilon}(t)&=-4\alpha\frac{1}{\epsilon}\int_{0}^{t}1\!_{\{0\leq\Delta X_{s}-2\alpha\Delta Z^{\alpha}_{s}\leq \epsilon\}}\left(1\!_{\{X^{1}_s \neq 0\}}-1\!_{\{X^{2}_s \neq 0\}}\right)\hspace{-.15cm}\left(1\!_{\{X^{1}_s > 0\}}-1\!_{\{X^{2}_s > 0\}}\right)ds.
\end{align*}
 Un calcul élémentaire nous donne
\begin{align*}
 \lim_{\epsilon\rightarrow 0}I^{1}_{\epsilon}(t)&=\lim_{\epsilon\rightarrow 0}\frac{1}{\epsilon}\left[\int_{0}^{t}\left(1\!_{\{X^{2}_{s}=0,\ 0< (1-2\alpha)X^{1}_s\leq \epsilon\}}+ 1\!_{\{X^{1}_{s}=0,\ -\epsilon\leq X^{2}_s< 0\}}\right)ds\right]\\
 &=\lim_{\epsilon\rightarrow 0}\left[\frac{1}{1-2\alpha}\frac{1}{\epsilon}\int_{0}^{t}1\!_{\{X^{2}_{s}=0,\ 0<X^{1}_s\leq \epsilon\}}ds+ \frac{1}{\epsilon}\int_{0}^{t}1\!_{\{X^{1}_{s}=0,\ -\epsilon\leq X^{2}_s< 0\}}ds\right]\\
 &=\frac{2\alpha}{1-2\alpha}\lim_{\epsilon\rightarrow 0}\frac{1}{\epsilon}\int_{0}^{t}1\!_{\{X^{2}_{s}=0,\ 0<X^{1}_s\leq\epsilon\}}ds +
 \lim_{\epsilon\rightarrow 0}\frac{1}{\epsilon}\int_{0}^{t}\left(1\!_{\{X^{2}_{s}=0,\ 0< X^{1}_s\leq \epsilon\}}+ 1\!_{\{X^{1}_{s}=0,\ -\epsilon\leq X^{2}_s< 0\}}\right)ds\\
 &=\frac{2\alpha}{1-2\alpha}\lim_{\epsilon\rightarrow 0}\frac{1}{\epsilon}\int_{0}^{t}1\!_{\{X^{2}_{s}=0,\ 0<X^{1}_s\leq \epsilon\}}ds +
 \lim_{\epsilon\rightarrow 0}\frac{1}{\epsilon}\int_{0}^{t}1\!_{\{0\leq\Delta X_{s}\leq \epsilon\}} \left[1\!_{\{X^{1}_s \neq 0\}}-1\!_{\{X^{2}_s \neq 0\}}\right]^{2}ds
\end{align*}
et comme $d\langle \Delta X,\Delta X \rangle  _{t}=\left[1\!_{\{X^{1}_t \neq 0\}}-1\!_{\{X^{2}_t \neq 0\}}\right]^{2}dt$, on en déduit

\begin{align}\label{terme1}
 \lim_{\epsilon\rightarrow 0}I^{1}_{\epsilon}(t)&=\frac{2\alpha}{1-2\alpha}\lim_{\epsilon\rightarrow 0}\frac{1}{\epsilon}\int_{0}^{t}1\!_{\{X^{2}_{s}=0,\ 0< X^{1}_s\leq \epsilon\}}ds +
 L_{t}^{0}(\Delta X).
 \end{align}
 De même
\begin{align}\label{terme2}
 \lim_{\epsilon\rightarrow 0}I^{2}_{\epsilon}(t)&=4\alpha^{2}\lim_{\epsilon\rightarrow 0}\frac{1}{\epsilon}\int_{0}^{t}1\!_{\{ 0\leq (1-2\alpha)X^{1}_s-X^{2}_{s}\leq \epsilon, X^{1}_s>0, X^{2}_s\leq0 \}}ds
 \end{align}
 et
 \begin{align}\label{terme3}
 \lim_{\epsilon\rightarrow 0}I^{3}_{\epsilon}(t)&=-4\alpha\lim_{\epsilon\rightarrow 0}\frac{1}{\epsilon}\int_{0}^{t}1\!_{\{ 0< (1-2\alpha)X^{1}_s\leq \epsilon, X^{2}_s=0 \}}ds
 \end{align}
D'où, en combinant (\ref{appro}), (\ref{terme1}), (\ref{terme2}) et (\ref{terme3}), on a
\begin{align*}
\int_{0}^{t}\varphi(s)dL_{s}^{0}(\Delta X)&=\frac{1}{2\alpha}[L_{t}^{0}(\Delta X)-L_{t}^{0}(\Delta Z^{\alpha})]\\
&=\frac{-1}{1-2\alpha}\lim_{\epsilon\rightarrow 0}\left(\frac{1}{\epsilon}\int_{0}^{t}1\!_{\{X^{2}_{s}=0,\ 0<X^{1}_s\leq \epsilon\}}ds\right)\\
&-2\alpha\lim_{\epsilon\rightarrow 0}\left(\frac{1}{\epsilon}\int_{0}^{t}1\!_{\{ 0\leq (1-2\alpha)X^{1}_s-X^{2}_{s}\leq \epsilon, X^{1}_s>0, X^{2}_s\leq0 \}}ds\right)\\
&+2\lim_{\epsilon\rightarrow 0}\left(\frac{1}{\epsilon}\int_{0}^{t}1\!_{\{0<X^{1}_s,\ 0<(1-2\alpha)X^{1}_s\leq \epsilon, X^{2}_s=0 \}}ds\right)
\end{align*}
ce qui fournit, en faisant tendre $\alpha$ vers zero
\[\int_{0}^{t}\varphi(s)dL_{s}^{0}(\Delta X)=\lim_{\epsilon\rightarrow 0}\frac{1}{\epsilon}\int_{0}^{t}1\!_{\{X^{2}_{s}=0,\ 0<X^{1}_s\leq \epsilon\}}ds\]
par ailleurs,  on a
\begin{align*}\lim_{\epsilon\rightarrow 0}\frac{1}{\epsilon}\int_{0}^{t}1\!_{\{X^{2}_{s}=0,\ 0<X^{1}_s\leq \epsilon\}}ds&=(1-2\alpha)\lim_{\epsilon\rightarrow 0}\frac{1}{\epsilon}\int_{0}^{t}1\!_{\{0<X^{1}_s,\ X^{2}_{s}=0,\ 0\leq(1-2\alpha)X^{1}_s\leq \epsilon\}}ds,
\end{align*}
qui tend vers zero quand $\alpha\rightarrow\frac{1}{2}$.
D'où
\[\int_{0}^{t}\varphi(s)dL_{s}^{0}(\Delta X)=0.\]
($\imath\imath$) Pour détérminer l'équation satisfaite par $h_{1,2}$, nous allons réutiliser les calcules entrepris pour
l'équation (\ref{int-diff}), mais cette fois pour le processus $X^{1}\vee X^{2}$. En effet,
on a d'une part,
\begin{eqnarray*}    L^{0}_{t}(X^{1}\vee X^{2})&=& \int_{0}^{t}1\!_{\{X^{1}_{s}\leq 0\}} dL^{0}_{s}(X^{2})+\int_{0}^{t}1\!_{\{X^{2}_{s}<0\}}dL^{0}_{s}(X^{1})\nonumber \\
&=&\int_{0}^{t}1\!_{\{X_{s}^{1}\vee X_{s}^{2}=0\}}\frac{ds}{\rho(h_s)}.
\end{eqnarray*}
et en utilisant le point ($\imath$) du lemme \ref{1}, on a
\begin{eqnarray*}
 \frac{L^{0}_{t}(X^{1}\vee X^{2})- L^{0-}_{t}(X^{1}\vee X^{2})}{2}&=&\dint_{0}^{t}1\!_{\{X_{s}^{1}\vee X_{s}^{2}=0\}}\alpha(s)\left[\!\!\!\phantom{\int}1\!_{\{X^{1}_{s}\leq 0\}}dL^{0}_{s}(X^{2})\right.\\
 &&\quad\left.+\frac{1}{2}1\!_{\{X^{2}_{s}<0\}}dL^{0}_{s}(X^{1})\right]
 +\frac{1}{2}\dint_{0}^{t}1\!_{\{X_{s}^{1}\vee X_{s}^{2}=0\}}dL_{s}^{0}(\Delta X)\nonumber\\
  &=&\dint_{0}^{t}\alpha(s)\left(1\!_{\{X^{1}_{s}\leq 0\}}dL^{0}_{s}(X^{2})+1\!_{\{X^{2}_{s}<0\}}dL^{0}_{s}(X^{1})\right)\\
  && \quad\quad\quad+\frac{1}{2}L_{t}^{0}(\Delta X)
\end{eqnarray*}
D'où, en tenant compte que $X^{1}$ et $X^{2}$ sont solutions de (\ref{W-G}), on obtient
\begin{eqnarray}\label{formule2}
 \nonumber
 \frac{L^{0}_{t}(X^{1}\vee X^{2})+ L^{0-}_{t}(X^{1}\vee X^{2})}{2}&=&L^{0}_{t}(X^{1}\vee X^{2})-\frac{L^{0}_{t}(X^{1}\vee X^{2})- L^{0-}_{t}(X^{1}\vee X^{2})}{2}\\
  &=&\int_{0}^{t}\frac{1-\alpha(s)}{\rho(h_s)}1\!_{\{X_{s}^{1}\vee X_{s}^{2}=0\}}ds-\frac{1}{2}L^{0}_{t}(\Delta X).
\end{eqnarray}
D'autre part, on montre facilement, comme dans la première partie, que la loi  \\
\noindent $\P_{x}(X_{t}^{1}\vee X_{t}^{2}\in dy)$ est absolument continue par rapport à
la mesure de Lebesgue sur $\R^{*}$ et que sa densité $\overline{p}_t(x,y)$ vérifie
\begin{eqnarray*}
  \frac{\overline{p}_{t}(x,0^{+})+ \overline{p}_{t}(x,0^{-})}{2}&=&\frac{1}{2}\frac{d}{dt}\E_{x}[L^{0}_{t}(X^{1}\vee X^{2})+L^{0-}_{t}(X^{1}\vee X^{2})]\\
  &=&\frac{1}{2\pi}\int \overline{f}(\lambda, t)d\lambda,
\end{eqnarray*}
où $\overline{f}(\lambda, t)=\E_{x}(e^{i\lambda X_t^{1}\vee X^{2}_t}1_{\{X_t^{1}\vee X^{2}_t\neq0\}})$. Et donc par (\ref{formule2}), on obtient
\begin{eqnarray}\label{form}
\frac{1-\alpha(t)}{\rho(h_t)}\!\ h_{1,2}(t) + \E_{x}\psi^{\prime}(t)=\frac{1}{2\pi}\int \overline{f}(\lambda, t)d\lambda
\end{eqnarray}
Par ailleurs, grâce à (\ref{supdecomp}), (\ref{bordsup}) et à la formule d'Itô on obtient
\begin{eqnarray*}
\overline{f}(\lambda,t)+ h_{1,2}(t)=e^{i\lambda x}&-&\frac{\lambda^{2}}{2}\dint_{0}^{t}(\overline{f}(\lambda,s)+h_{1,2}(s))ds
+\dint_{0}^{t}\left[\frac{i\lambda}{\rho(h_s)}\alpha(s)+ \frac{\lambda^{2}}{2}\right]h_{1,2}(s) ds\\
&+& \frac{i\lambda}{2}\E_{x}[L^{0}_{t}(\Delta X)]
\end{eqnarray*}
ce qui donne
\begin{eqnarray}
\overline{f}(\lambda,t)= e^{-\frac{\lambda^{2}}{2}t}\left[e^{i\lambda x}+\overline{u}(t)\right]- h_{1,2}(t),
\end{eqnarray}
avec
 \[\overline{u}(t)=\dint_{0}^{t}ds\ e^{\frac{\lambda^{2}}{2}s}\left[\frac{i\lambda}{\rho(h_s)}\alpha(s)h_{1,2}(s)+ \frac{\lambda^{2}}{2}h_{1,2}(s)-i\lambda\psi^{\prime}(s)\right].\]

\noindent Et par suite, en utilisant le théorème de Fubini et le fait que la fonction
\[\lambda\mapsto \frac{i\lambda}{\rho(h_{t-s})}e^{\frac{\lambda^{2}}{2}s}\alpha(t-s)h_{1,2}(t-s)-i \lambda\psi^{\prime}(s)\]
est impaire, on a
\begin{eqnarray*}
\frac{1}{2\pi}\int \overline{f}(\lambda, t)d\lambda&=&\frac{e^{-\frac{x^{2}}{2 t}}}{\sqrt{2\pi t}}+\frac{1}{2\pi}\int d\lambda\left[\int_{0}^{t}dse^{-\frac{\lambda^{2}}{2}(t-s)}\left(\frac{i\lambda}{\rho(h_s)}\alpha(s)h_{1,2}(s)+\frac{\lambda^{2}}{2}h_{1,2}(s)\right.\right.\\
&&\quad\left.\left.-i \lambda\psi^{\prime}(s)\!\!\!\!\!\!\!\phantom{\int}\right)-h_{1,2}(t)\!\!\!\!\!\phantom{\int}\right]\\
&=&\frac{e^{-\frac{x^{2}}{2 t}}}{\sqrt{2\pi t}}+\frac{1}{2\pi}\int d\lambda\left[\int_{0}^{t}dse^{-\frac{\lambda^{2}}{2}(t-s)}\frac{\lambda^{2}}{2}h_{1,2}(s)-h_{1,2}(t)\right]\\
&=&\frac{1}{\sqrt{2\pi t}}\!\ e^{-\frac{x^{2}}{2 t}}+\frac{1}{2\sqrt{2\pi}}\dint_{0}^{t} \frac{h_{1,2}(t-s)-h_{1,2}(t)}{s^{3/2}}\! ds -\frac{h_{1,2}(t)}{\sqrt{2\pi t}},
\end{eqnarray*}
ce qui joint à (\ref{form}) prouve (\ref{int-diff2}).\\
Cherchons maintenant l'équation satisfaite par $\overline{h}_{1,2}$. Pour cela, écrivons $\overline{h}_{1,2}(t)=\P_{x}(Y_{t}^{1}\vee Y_{t}^{2}=0)$, avec $Y^{i}=-X^{i}$, et observons que puisque
\[L^{0}_{t}(X^{i})=L^{0-}_{t}(-X^{i})\quad \mbox{et}\quad L^{0}_{t}(-X^{i})-L^{0-}_{t}(-X^{i})=2\int_{0}^{t}\alpha(s)dL^{0}_{s}(X^{i})\]
alors
\[dL^{0}_{t}(X^{i})= \frac{1}{1-2\alpha(s)}dL^{0}_{t}(Y^{i})\]
ce qui implique que les $Y^{i}$, $i=1,2$, satisfont aux équations
\begin{align}\label{Y}
 \left\{
\begin{aligned}
&dY^{i}_{t} =1\!_{\{Y^{i}_t \neq 0\}}d\beta_t + \widetilde{\alpha}(t) dL^{0}_t(Y^{i})\\
&\int_{0}^{t}1\!_{\{Y^{i}_s= 0\}}ds =\int_{0}^{t}\varrho(s)dL^{0}_{s}(Y^{i}),
\end{aligned}
\right.
\end{align}
avec $\widetilde{\alpha}(t)=-\frac{\alpha(t)}{1-2\alpha(t)}$, $\varrho(t)=\frac{\rho(h_t)}{1-2\alpha(t)}$ et  $\beta=-B$. Pour conclure à la démonstrations, il suffit alors de remarquer que
\begin{eqnarray*}
 \nonumber
 \frac{L^{0}_{t}(Y^{1}\vee Y^{2})+ L^{0-}_{t}(Y^{1}\vee Y^{2})}{2}&=&\int_{0}^{t}\frac{1-\widetilde{\alpha}(s)}{\varrho(s)}1\!_{\{Y_{s}^{1}\vee Y_{s}^{2}=0\}}ds-\frac{1}{2}L^{0}_{t}(X^{2}-X^{1})\\
 &=&\int_{0}^{t}\frac{1-\alpha(s)}{\rho(h_s)}1\!_{\{Y_{s}^{1}\vee Y_{s}^{2}=0\}}ds-\frac{1}{2}L^{0}_{t}(X^{2}-X^{1}).
\end{eqnarray*}
et que $L^{0}_{t}(X^{2}-X^{1})=L^{0-}_{t}(\Delta X)=L^{0}_{t}(\Delta X)$ puisque la partie à variation finie de $X^{1}-X^{2}$ ne charge pas l'ensemble $\{t, X_t^{1}=X_t^{2}\}$.

 \end{proof1}

\noindent \begin{proof2}[Preuve de la Proposition \ref{tempsup}] L'expression du temps local $L^{0}_{t}(X^{1}\vee X^{2})$ résult directement de (\ref{templocalsup1}) et (\ref{intphi}). Montrons maintenant que  $L_{t}^{0}(\Delta X)=0$.\\
\noindent Observons que, pour tout $t\geq0$, $$2h(t)=h_{1,2}(t)+\overline{h}_{1,2}(t).$$
On en déduit alors, en combinant (\ref{int-diff2}) et (\ref{int-diff3}), que $h$ est solution de
\begin{eqnarray}\label{int-diff4}
  \frac{1-\alpha(t)}{\rho(h_t)}\!\ h(t)+ \E\psi^{\prime}(t)=\frac{1}{\sqrt{2\pi t}}\!\ e^{-\frac{x^{2}}{2 t}} +\hspace{-.1cm} \frac{1}{2\sqrt{2\pi}}\hspace{-.08cm} \dint_{0}^{t} \frac{h(t-s)-h(t)}{s^{3/2}}\!\ ds -\hspace{-.08cm} \frac{h(t)}{\sqrt{2\pi t}}
\end{eqnarray}
ce qui entraîne, grâce à (\ref{int-diff}), que
\[0=\E\psi^{\prime}(t)=-\frac{1}{2}\frac{d}{dt}\E L_{t}^{0}(\Delta X)\]
et donc, $L_{t}^{0}(\Delta X)=0$. Ce qui achève la démonstration.
\end{proof2}


\end{document}